\documentclass[conference]{IEEEtran}
\IEEEoverridecommandlockouts
\usepackage{cite}
\usepackage{amsmath,amssymb,amsfonts}
\usepackage{graphicx}
\usepackage{textcomp}
\usepackage{xcolor}

\usepackage{dsfont} 
\usepackage{url}
\usepackage[utf8]{inputenc} 
\usepackage[format=hang,singlelinecheck=false,caption=false,font=footnotesize]{subfig}
\usepackage{changes}
\usepackage{algpseudocode}

\def\mathrlapinternal#1#2{\rlap{$\mathsurround=0pt#1{#2}$}}
\def\mathrlap{\mathpalette\mathrlapinternal}

\newcommand{\R}{\mathbb{R}}

\renewcommand{\epsilon}{\varepsilon}
\renewcommand{\phi}{\varphi}

\newcommand{\reno}{\mathit{Re}}
\newcommand{\derivoper}{\mathrm{d}}
\newcommand{\fderivoper}{\mathcal{D}}
\newcommand{\one}{\mathds{1}}

\DeclareMathOperator{\pmv}{pmv}

\newcommand{\abs}[1]{\left\lvert #1\right\rvert}
\newcommand{\absb}[1]{\bigl\lvert #1\bigr\rvert}

\newcommand{\norm}[1]{\left\lVert #1\right\rVert}
\newcommand{\normb}[1]{\bigl\lVert #1\bigr\rVert}

\newcommand{\set}[1]{\left\{ #1\right\}}
\newcommand{\setb}[1]{\bigl\{ #1\bigr\}}

\newcommand{\setbb}[1]{\biggl\{ #1\biggr\}}

\newcommand{\paren}[1]{\left( #1\right)}
\newcommand{\parenb}[1]{\bigl( #1\bigr)}
\newcommand{\parenB}[1]{\Bigl( #1\Bigr)}
\newcommand{\parenbb}[1]{\biggl( #1\biggr)}
\newcommand{\parenBB}[1]{\Biggl( #1\Biggr)}
\newcommand{\parenn}[1]{( #1)}

\let\p=\paren
\let\pn=\parenn
\let\pb=\parenb
\let\pB=\parenB

\newcommand{\sqparen}[1]{\left[ #1\right]}
\newcommand{\sqparenb}[1]{\bigl[ #1\bigr]}
\newcommand{\sqparenB}[1]{\Bigl[ #1\Bigr]}
\newcommand{\sqparenbb}[1]{\biggl[ #1\biggr]}
\newcommand{\sqparenBB}[1]{\Biggl[ #1\Biggr]}
\newcommand{\sqparenn}[1]{[ #1]}

\let\spb=\sqparenb

\newcommand{\dprod}[2]{\left\langle #1, #2\right\rangle}
\newcommand{\dprodb}[2]{\bigl\langle #1, #2\bigr\rangle}

\let\smat=\sqmatrix
\let\smats=\sqmatrixs

\newcommand{\deriv}[3][]{\frac{\derivoper^{#1} #2}{\derivoper #3^{#1}}}
\newcommand{\pderiv}[3][]{\frac{\partial^{#1} #2}{\partial #3^{#1}}}
\newcommand{\fderiv}[2][]{\fderivoper_{#1} #2}
\newcommand{\diff}[1]{\derivoper #1}

\newcommand{\diver}[2][]{\nabla_{\!#1} \cdot #2}

\newcommand{\grad}[2][]{\nabla_{\!#1} #2}
\newcommand{\lapl}[2][]{\triangle_{#1} #2}
\newcommand{\bdry}[1]{\partial #1}
\newcommand{\unit}[1]{\,{\rm #1}}
\newcommand{\unitns}[1]{\ensuremath{{\rm #1}}}

\newcommand{\uf}[1]{\textbf{DEPRECATED}}

\newtheorem{theorem}{Theorem}
\newtheorem{lemma}{Lemma}
\newtheorem{definition}{Definition}

\def\BibTeX{{\rm B\kern-.05em{\sc i\kern-.025em b}\kern-.08em
    T\kern-.1667em\lower.7ex\hbox{E}\kern-.125emX}}
\begin{document}

\title{Comfort-Aware Building Climate Control Using Distributed-Parameter Models
}

\author{\IEEEauthorblockN{Runxin He}
\IEEEauthorblockA{\textit{Department of Electrical \&
Systems Engineering} \\
\textit{Washington University in St. Louis}\\
Saint Louis, USA \\
runxinhe@hotmail.com}
\and
\IEEEauthorblockN{Humberto Gonzalez}
\IEEEauthorblockA{\textit{LinkedIn} \\ 
Sunnyvale, USA \\
hgonzale@wustl.edu}}

\maketitle
\begin{abstract}
Controlling Heating, Ventilation and Air Conditioning (HVAC) system to maintain occupant's indoor thermal comfort is important to energy-efficient buldings and the development of smart cities.
In this paper, we formulate a model predictive controller (MPC) to make optimal control strategies to HVAC in order to maintain occupant's comfort by predicted mean vote index, and then present a whole system with an estimator to indoor climate and apartment's geometric information based on only thermostats.
In order to have accurate spatial resolution and make the HVAC system focus on only a zoned area around the occupant, a convection-diffusion Computer Fluid Dynamics (CFD) model is used to describe the indoor air flow and temperature distribution.
The MPC system generates corresponding PDE-contrained optimization problems, and we solve them by obtain the gradients of cost functions with respect to problems' variables with the help of CFD model's adjoint equations.
We evaluate the performance of our method using simulations of a real apartment in the St.\ Louis area.
Our results show our MPC system's energy efficiency and the potential for its application in real-time operation of high-performance buildings.
\end{abstract}

\begin{IEEEkeywords}
Model Predictive Control,
Computational Fluid Dynamics,
PDE Optimization,
Thermal Comfort
\end{IEEEkeywords}

\section{Introduction}
\label{sec:intro}
%
%
%
%

Heating, Ventilation, and Air Conditioning (HVAC) units are complex systems that coordinate incorporate mechanical and electrical components using data processing algorithms, designed to control the climate within buildings.
Standard HVAC systems typically use a small number of thermostats, together with simple temperature regularization loops~\cite{freire2008predictive}, to follow a temperature set point and ignore other variables, such as humidity and flow velocity.
The multi-dimensional nature of human comfort lead to the adoption of the Predicted Mean Vote (PMV) index~\cite{fanger1970thermal} in the the ASHRAE~55 Standard~\cite{ashrae2010} to quantify indoor thermal comfort.
Our effort in this paper is motivated by our interest to explicitly incorporate PMV-based comfort measurements to the HVAC actuation loop with distributed-parameter Computer Fluid Dynamic (CFD) model.

Most studies in HVAC control can be categorized within two main trends.
The first corresponds to learning-based methods~\cite{yuce2014utilizing, huang2015neural}, which simplify the dynamic complicity by large quantities of data and training time.
The second trend corresponds to model-based predictive methods, which describe the dynamic evolution of building climate variables using first-principle physical models, or their approximations.
Among model-based predictive methods, MPC stands out to be one of the standard methods due to its flexible mathematical formulation and its robust performance in real-world implementations~\cite{huang2011model, xi2007support, morocsan2010building, afram2014theory}.
MPC has been applied to zoned temperature control and temperature regularization in the past~\cite{ma2012demand, chandan2013optimal, he2015zoned}, showing significant improvements in energy efficiency compared to other classical control methods~\cite{privara2011model, sturzenegger2016model}.
Also, the flexibility of MPC allows for the stochastic nature of certain disturbance variables, such as weather forecast~\cite{oldewurtel2012use}.

There exist many results regarding the use of indoor climate CFD models in optimal control schemes.
For example, Doering and Gibbon~\cite{doering1995applied} and Ito~\cite{ito1998optimal} studied the existence and smoothness of fluid models under certain conditions.
In the specific case of HVAC control, Burns et al.~\cite{borggaard2009control, burns2012control, burns2013approximation} used linearized CFD models and quadratic cost functions.
Their results show that these PDE models 
can be effectively used for building control, while the linearized approximations are accurate only if it is close to the original steady-state linearization point~\cite{he2015zoned}.

In this paper we focus on efficiently utilizing existing HVAC system technologies to maximize user comfort, rather than simply trying to maintain a constant temperature. To achieve our goal we develop estimation and control algorithms that consider spatio-temporal distributions of temperature and air flow, which take into account changes in floor-plan geometry (such as doors being opened or closed), outdoor weather, and the position of fans and portable heaters.
All this information is fed into a distributed-parameter nonlinear CFD model and a multi-dimensional human comfort index to generate, and optimize, short-term predictions.
Following our results in~\cite{he2016gradient}, we develop a first-order gradient-based optimization method to find optimal solutions using these nonlinear CFD models.
Our paper is organized as follows.
Sec.~\ref{sec:prb_state} describes the CFD model and the formulation of our MPC scheme.
Sec.~\ref{sec:methods} presents our gradient-based algorithm to solve the PDE-constrained optimization problem, as well as our finite-element discretization used numerical computations.
Sec.~\ref{sec:results} presents the results of our simulated experiments, based on a real apartment in the St.\,Louis area.

\section{PROBLEM STATEMENT}
\label{sec:prb_state}
\begin{figure}[tp]
  \centering
  \includegraphics[width=.8\linewidth]{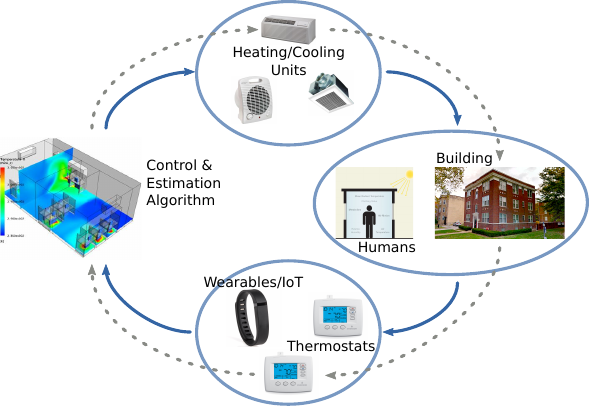}
  \caption{%
    Climate control solutions.
    Traditional solutions (shown in dashed gray lines) involve a single thermostat, one ventilation unit, and no human comfort input.
    Our new paradigm (shown in solid blue) considers multiple heterogeneous sensors, a controller and estimator based on distributed-parameter models, and multiple ventilation units.
  }
  \label{fig:mpc_comfort}
\end{figure}

We consider the problem of controlling the indoor temperature and air flow in a building equipped with a HVAC system consisting of one or several small units, with the goal of simultaneously improving both energy efficiency and human comfort.
As illustrated in Fig.~\ref{fig:mpc_comfort}, traditional solutions to this problem rely on single thermostat measurements and a single-unit HVAC systems, thus equating human comfort to a single temperature value in the building as measured by the thermostat.
Instead, we propose to use heterogeneous Internet-of-Things sensor data~\cite{jin2014environment, cleland2013optimal, yang2010review} together with first-principle fluid dynamics models to estimate and predict the spatial and temporal behavior of all the relevant climate variables, which in turn allows us to achieve our goals of efficiency and comfort.
Our improved climate control scheme is also illustrated in Fig.~\ref{fig:mpc_comfort}.

In this section we describe in detail all the components in our control scheme, including the CFD model, comfort evaluation model, and optimization algorithms for control and estimation.

\subsection{Computational Fluid Dynamic Model}
\label{subsec:cfd_model}
Before formally presenting our control scheme, we define the notation used throughout this paper.
We will endow $\R^n$ with the traditional 2-norm, denoted $\norm{\cdot}$, and we will denote the inner product between $x_1$ and $x_2$ in $\R^n$ by $x_1 \cdot x_2$.

Our functional analysis notation closely follows that in~\cite{adams2003sobolev}.
Let $S$ be a subset of a vector space, then we say that $L^2(S)$ is the space of squared-integrable functions with domain in $S$, endowed with the inner product $\dprod{f_1}{f_2}_S = \int_{S} f_1(x) \cdot f_2(x)\, \diff{x}$, and the norm $\norm{f}_{S} = \sqrt{\dprod{f}{f}}$.
Moreover, $H^1(S)$ is the Sobolev space of functions in $L^2(S)$ whose weak-derivatives are also in $L^2(S)$, and $H_0^1(S)$ is the set of functions $f \in H^1(S)$ whose boundary is trivial, i.e., $f(x) = 0$ for each $x \in \bdry{S}$.
We define $L^2\pb{[t_0, t_f]; H^1(S)}$ as the set of functions $u\colon [t_0,t_f] \to H_0^1(S)$ such that $\int_{t_0}^{t_f} \norm{u(t)}_S^2 + \normb{\deriv{u(t)}{x}}_S^2\, \diff{t}$ is finite.

The kernel of our CFD model is the incompressible Navier-Stokes equation, which is a good approximation for the coupling of temperature with free flow convection at atmospheric conditions~\cite{awbi1989, Dobrzynski2004}.
Throughout the paper we make two major simplifications to fluid dynamic model.
First, we assume that the air flow behaves as a laminar fluid which reaches steady state much faster than the temperature.
Both theoretical~\cite{awbi2003ventilation} and experimental~\cite{sun2003} results have shown that turbulent flows are present in residential building, yet their overall effects are negligible.
Second, we consider only two-dimensional air flows moving parallel to the ground.
Both assumptions reduce the accuracy of our model to some extent~\cite{van2013comparison}, yet they allow us to significantly simplify the computational complexity.

Following the notation presented above, let $\Omega \subset \R^2$ be the region of interest, assumed to be compact and connected.
Let $u\colon \Omega \to \R^2$ be the \emph{stationary air flow velocity}, $p\colon \Omega \to \R$ be the \emph{stationary air pressure}, and $T_e\colon \Omega \times [t_0, t_f] \to \R$ be the \emph{temperature}.
Then, the temperature convection-diffusion model~\cite{landau_1971} over the region $\Omega$ can be described by the following PDE:
\begin{multline}
  \label{eq:temp}
  \pderiv{T_e}{t}(x,t)
  - \grad[x] \cdot \pb{\kappa(x)\, \grad[x]{T_e}(x,t)}
  + u(x) \cdot \grad[x]{T_e}(x,t) =\\
  = g_{T_e}(x,t),
\end{multline}
where $g_{T_e}\colon \Omega \times [t_0, t_f] \to \R$ represents the heat source in the room, $\kappa\colon \Omega \to \R$ is the \emph{thermal diffusivity}, and $\grad[x]{} = \pb{\pderiv{}{x_1},\pderiv{}{x_2}}^T$ is the \emph{gradient operator}.
The initial condition of the temperature is:
\begin{equation}
  \label{eq:init_t}
  T_e(x, t_0) = \pi_0(x),\quad \forall x \in \Omega.
\end{equation}
Similarly, the air flow in $\Omega$ is governed by the following stationary incompressible Navier-Stokes PDEs:
\begin{align}
  \notag
  - \frac{1}{\reno} \lapl[x]{u}(x)
  + \pb{u(x) \cdot \grad[x]{}}\, u(x)
  + \grad[x]{p}(x) &+{} \\
  \label{eq:ns1}
  {}+ \alpha(x)\, u(x) &= g_u(x), \\
  \label{eq:ns2}
  \diver[x]{u}(x) &= 0,
\end{align}
where $g_u \colon \Omega \to \R^2$ represents the external forces applied to the air (such as fans), $\reno$ is the \emph{Reynolds number} of the air, and $\lapl[x]{} = \pderiv[2]{}{x_1} + \pderiv[2]{}{x_2}$ is the \emph{Laplacian operator}.
We introduce the \emph{viscous friction coefficient} $\alpha\colon \Omega \to \R$~\cite{gersborg2005topology}, to model the effect of solid materials in the air flow.
Indeed, when the point $x$ corresponds to a material that blocks air, we choose $\alpha(x) \gg u(x)$, which results in $u(x) \approx 0$, and when the point $x$ corresponds to air we choose $\alpha(x) = 0$.

We assume that the exterior walls of the building are solid and correspond to the boundary of $\Omega$, denoted $\bdry{\Omega}$.
Thus, the boundary condition for the temperature is:
\begin{equation}
  \label{eq:bnd_t}
  T_e(x) \equiv T_A,\quad \forall x \in \bdry{\Omega},
\end{equation}
where $T_A$ is the \emph{atmospheric temperature}, and the boundary condition for the air flow is:
\begin{equation}
  \label{eq:bnd_u1}
  u(x) \equiv 0,\quad \forall x \in \bdry{\Omega}.
\end{equation}
Moreover, $p(x) \equiv p_A$ for each $x \in \bdry{\Omega}$, where $p_A$ is the \emph{atmospheric pressure}.

The existence and uniqueness of weak solution to the CFD model in equations~\eqref{eq:temp} to~\eqref{eq:bnd_u1} are given formalized below.

\begin{theorem}[Existence of weak solutions]
  \label{thm:wk_existence1}
  Let $\alpha, \kappa \in L^2(\Omega)$ be bounded functions, $g_{T_e} \in L^2\pb{[t_0,t_f]; H^1(\Omega)}$, and $g_u \in L^2(\Omega) \times L^2(\Omega)$.
  Then, there exists weak solutions $T_e \in L^2\pb{[t_0,t_f]; H^1(\Omega)}$, $u \in L^2(\Omega) \times L^2(\Omega)$ and $p \in L^2(\Omega)$ to the PDEs in eqs.~\eqref{eq:temp} to~\eqref{eq:bnd_u1}.
\end{theorem}

\begin{theorem}[Uniqueness of weak solutions]
  \label{thm:wk_unique1}
  Let $\abs{\Omega}$ denote the area of $\Omega$, and consider the same conditions as in Thm.~\ref{thm:wk_existence1}.
  If $\norm{\grad[x]{u}}_{L^2} < \frac{2}{\abs{\Omega}^{1/2}\, \reno}$, then the weak solutions $T_e$, $u$, and $p$, defined in Thm.~\ref{thm:wk_existence1}, are unique.
\end{theorem}

The proofs of Thms.~\ref{thm:wk_existence1} and~\ref{thm:wk_unique1}
are omitted from this section, and presented instead in full detail in Appendix~\ref{appendix:existence}.
Note that the condition in Thm.~\ref{thm:wk_unique1} is simply a sufficient condition, which corresponds to the state-of-the-art in terms of uniqueness conditions for stationary Navier-Stokes equations.

\medskip
We assume that there are $n_t$ temperature sensors in the building to estimate both the indoor fluid dynamic and floor configuration change.
The detailed model for the sensors and related estimation problem are discussed in~\cite{he2016gradient}.
We also consider $n_f$ fan units distributed in different locations of the building, capable of blowing air while introducing or removing heat from the environment.
Furthermore, we assume that there are $n_d$ doors in the building.
We define $\theta_i \in \set{0,1}$ as the configuration of the $i$-th door, i.e., $\theta_i = 1$ when the $i$-th door is open, and $\theta_i = 0$ when it is closed.
Let $\Omega_{\theta_i} \subset \Omega$ be the area occupied by the $i$-th door when it is closed.
When the door configuration changes by the estimation problem~\cite{he2016gradient}, so does the prediction generated by our CFD model.
In particular, the parameters $\alpha$ and $\kappa$ change for each $x \in \Omega_{\theta_i}$ as a function of $\theta_i$.
We model this relation by defining $\alpha\colon \Omega \times \set{0,1}^{n_d} \to \R$ and $\kappa\colon \Omega \times \set{0,1}^{n_d} \to \R$ as follows:
\begin{equation}
  \label{eq:kappa}
  \begin{aligned}
  \alpha(x,\theta)
  &= \alpha_w\, \sum_{i=1}^{n_d} \p{1-\theta_i}\, \one\set{x \in \Omega_{\theta_i}},
  \quad \text{and} \\
  \kappa(x,\theta)
  &= \kappa_0 + \p{\kappa_w-\kappa_0}\, \sum_{i=1}^{n_d} \p{1-\theta_i}\, \one\set{x \in \Omega_{\theta_i}},
  \end{aligned}
\end{equation}
where $\kappa_0$ is the open air thermal diffusivity, $\alpha_w$ and $\kappa_w$ are the solid wall viscosity and diffusivity, and $\one\set{\cdot}$ is the indicator function.
Note that both $\alpha$ and $\kappa$ are affine functions of $\theta \in \R^{n_d}$.

\subsection{Predicted Mean Vote Index Approximation}
\label{subsec:pmv}
The PMV index was proposed by Fanger~\cite{fanger1970thermal} and recommended by ASHRAE~\cite{ashrae2010} in order to predict the average vote of a large group of persons on the thermal sensation scale.
It uses heat balance equation to relate six key factors to the average response of people on the thermal comfort.
According to~\cite{ashrae2010}, the closer the PMV index to zero the more comfortable the occupant feels; and the range of PMV index for an acceptable thermal environment of general comfort is from $-0.5$ to $0.5$.


For each person in the building, the PMV index is an implicit function
of the following six factors: metabolic rate, clothing insulation, air
temperature and humidity, air velocity, and mean radiant temperature.
Due to its implicit formulation, and in order to reduce the
computational effort in our optimization algorithms, we apply known
approximations of the PMV index~\cite{cigler2012optimization,
  federspiel1992}. 
A detailed presentation of the PMV index approximations is beyond the
scope of this section, but it is presented in detail in Appendix~\ref{appendix:pmv}.
After applying these approximations, the PMV index is influenced by
both the resident's personal information and variables in our CFD
model.
Given personal information, define the index value as
\begin{equation}
  \begin{aligned}
    & pmv(T_e(x,t;g_{T_e},g_u), u(x;g_u))
    = a_0
    + \frac{a_1 g_u + a_2}{a_3 g_u} + \\
    & \hspace{1.0em}
    + a_4 T_e(x,t;g_{T_e},g_u)
    + a_5 T_e(x,t;g_{T_e},g_u)^2 + \\
    & \hspace{1.0em}
    + \frac{a_6 + a_7 T_e(x,t;g_{T_e},g_u)}
    {a_8 + a_{9} T_e(x,t;g_{T_e},g_u) + a_{10} T_e(x,t;g_{T_e},g_u)^2}.
  \end{aligned}
  \label{eq:pmv_simplification}
\end{equation}
Where $\{a_i\}_{i=0}^{10}$ are parameters related to residents's
personal information, they are presented in detail in Appendix~\ref{appendix:pmv}.

\subsection{Model Predictive Control and Estimation}
\label{subsec:optimization}
Model predictive control and estimation is a scheme where two optimization problems, one to estimate the initial condition of the system and its unknown parameters, the second to solve a finite-horizon optimal control, are iteratively applied using a receding horizon approach.
Given a sequence of times, $\set{t_i}_{i=0}^\infty$, at each $t_i$ we solve an optimal estimation problem over the past horizon $[t_i - T, t_i]$, and an optimal control problem over the future horizon $[t_i, t_i + T']$.
Note that, in practice, the sequence $\set{t_i}_{i=0}^\infty$ is typically equidistant with $t_i - t_{i-1} \ll \min\set{T,T'}$.
The solution of the optimal control problem solved at $t_i$ is applied in the interval $[t_i,t_{i+1}]$.

First, we formulate our optimal estimation problem to find the door configuration $\theta$ and initial temperature $\pi_0$, using the information from the $n_t$ thermostats in the building.
Given the estimation time horizon as $[t_i - T, t_i]$, $\pi_0 \colon \Omega \to \R$, and $\theta \in \R^{n_d}$, 
we define our optimal estimation problem as
\begin{equation}
  \label{eq:cost_estimate}
  \begin{aligned}
    \min_{\pi_0,\, \theta}\quad
    & J_e\p{T_e, u, p, \pi_0, \theta},\\
    \text{subject to:}\ \
    & \text{PDE constraints~\eqref{eq:temp}, \eqref{eq:ns1}, and~\eqref{eq:ns2},} \\
    & \text{initial constraints~\eqref{eq:init_t},} \\
    & \text{boundary constraints~\eqref{eq:bnd_t}, and~\eqref{eq:bnd_u1},} \\
    & \theta_i \in [0,1], \quad \forall i \in \set{1,\dotsc,n_d}.
  \end{aligned}
\end{equation}
For further reading, the detailed descriptions to the problem~\eqref{eq:cost_estimate} is in~\cite{he2016gradient}.

Second, we formulate the optimal control problem to find the fan and heater control signals, denoted $g_{T_e}$ and $g_u$ in eqs.~\eqref{eq:temp} and~\eqref{eq:ns1}.
The optimal control problems takes the optimal estimations computed by the problem in eq.~\eqref{eq:cost_estimate}, denoted $\hat{\pi}_0$ and $\hat{\theta}$.
Abusing notation, let $\widehat{T}_e(x) = T_e(x, t_i; \hat{\pi}_0, \hat{\theta})$, i.e., the temperature distribution at time $t_i$ obtained from the optimal result of the problem in eq.~\eqref{eq:cost_estimate}.
Given the optimal control horizon as $[t_i, t_i + T']$, $g_{T_e}\colon [t_i, t_i + T'] \times \Omega \to \R$, and
$g_u\colon \Omega \to \R^2$, we define our optimal control objective as
\begin{multline}
  \label{eq:mpc_objective}
  \hspace{-1em}
    J_c
    \!=\!\! \int_{\Omega_t} \!\int_{t_i}^{\mathrlap{t_i + T'}} \abs{\pmv\pb{T_e(x,t;g_{T_e},g_u), u(x;g_u)}}^2\! \diff{t}\, \diff{x} +\\
    + \eta_0'\, \norm{g_{T_e}}_{\Omega \times [t_i,t_i+T']}^2
    + \eta_1'\, \norm{g_u}_{\Omega}^2,
\end{multline}
where $\eta'_0, \eta_1' > 0$ are weight parameters, $\Omega_t \subset \Omega$ is the target area of interest within the building, and $T_e(x,t;g_{T_e},g_u)$ and $u(x;g_u)$ are the solutions of eqs.~\eqref{eq:temp}, \eqref{eq:ns1}, and~\eqref{eq:ns2} with control signals $g_{T_e}$ and $g_u$.
Note that while in most applications the target area $\Omega_t = \Omega$, our formulation enables us to consider the case of smart zoned climatization, for example by only considering the rooms currently being occupied.
Using the objective function in eq.~\eqref{eq:mpc_objective} and the results of the optimal estimator, $\widehat{T}_e$ and $\hat{\theta}$, we formulate our optimal control problem as
\begin{equation}
  \begin{aligned}
    \label{eq:cost_control}
    \min_{g_{T_e}, \, g_u}\quad
    & J_c \p{T_e, u, g_{T_e}, g_u}, \\
    \text{subject to:}\ \
    & \text{PDE constraints~\eqref{eq:temp}, \eqref{eq:ns1}, and~\eqref{eq:ns2},} \\
    & \text{boundary constraints~\eqref{eq:bnd_t} and~\eqref{eq:bnd_u1},} \\
    & \pi_0(x) = \widehat{T}_e(x),\ \forall x \in \Omega, \\
    & g_{T_e}(x,t) \in \spb{\underline{g_{T_e}}, \overline{g_{T_e}}},\ \forall x \in \Omega,\ t \in [t_i,t_i+T'], \\
    & g_u(x) \in \spb{\underline{g_u}, \overline{g_u}},\ \forall x \in \Omega,
  \end{aligned}
\end{equation}
where $\underline{g_{T_e}}$, $\overline{g_{T_e}}$, $\underline{g_u}$ and $\overline{g_u}$ are minimum and maximum heater and fan power parameters, respectively.

\section{NUMERICAL METHOD}
\label{sec:methods}
Since the numerical algorithm to solve problem~\eqref{eq:cost_estimate} has been discussed in~\cite{he2016gradient}, in this section we develop a numerical algorithm to solve the optimization problem~\eqref{eq:cost_control} defined in Sec.~\ref{subsec:optimization} by adjoint equations, then present the whole system design to implement model predictive HVAC estimation and control.

\subsection{Adjoint Variables and Fr\'echet Derivatives}
\label{sec:adj_eq}

Let $\set{\lambda_i}_{i=1}^6$ be the adjoint variables, or Lagrange multipliers, each associated to one of the eqs.~\eqref{eq:temp} to~\eqref{eq:bnd_u1} and defined in its respective dual space.
Then, the Lagrangian function~\cite{giles1997adjoint, gunzburger2000adjoint} of our optimal estimation problem is:
\begin{multline}
  \label{eq:lg}
  L\pb{T_e, u, p, \pi_0, \theta, g_{T_e}, g_u, \set{\lambda_i}_{i=1}^6} ={} \\
    = \dprodb{\lambda_1}{\pderiv{T_e}{t} - \grad[x] \cdot \p{\kappa\, \grad[x]{T_e}} + u \cdot \grad[x]{T_e} - g_{T_e}}_{\Omega \times [t_0, t_f]} +{} \\
    + \dprodb{\lambda_2}{- \frac{1}{\reno} \lapl[x]{u} + \p{u \cdot \grad[x]{}} u + \grad[x]{p} + \alpha\, u - g_u}_{\Omega} +{} \\
  + \dprod{\lambda_3}{\diver[x]{u}}_{\Omega}
  + \dprod{\lambda_4}{T_e}_{\bdry{\Omega} \times [t_0, t_f]}
  + \dprod{\lambda_5}{u}_{\bdry{\Omega}} +{} \\
  + \dprod{\lambda_6}{T_e(\cdot,t_0) - \pi_0}_{\Omega}
  + J_c\p{g_{T_e}, g_u}.
\end{multline}
A necessary condition for optimality is that the inner product of the partial derivatives of $L$ with respect to each primal variable is equal to zero~\cite{gunzburger2000adjoint}.
After manipulating the resulting partial derivatives, this necessary condition can be transformed into a new set of PDEs that the adjoint variables must satisfy.
\begin{theorem}
  \label{thm:estimate_lg1}
  Let $\p{T_e^*, u^*, p^*, g_{T_e}^*, g_u^*}$ be a minimizer of $J_c$ subject to eqs.~\eqref{eq:temp} to~\eqref{eq:bnd_u1}.
 Then, there exist adjoint variables, $\set{\lambda_i}_{i=1}^6$, defined in their corresponding dual spaces as shown in eq.~\eqref{eq:lg}, which are weak solutions of the following equations for $x \in \Omega$ and $t \in [t_i, t_i + T']$:
  \begin{align}
    \label{eq:ad_pde1}
      - \grad[T_e]{J_c}
      + \pderiv{\lambda_1}{t}
      + \grad[x]{} \cdot \p{\kappa\, \grad[x]{\lambda_1}}
      + u^* \cdot \grad[x]{\lambda_1}
      &= 0, \\
    \notag
    \int_{t_i}^{\mathrlap{t_i + T'}} \lambda_1\, \grad[x]{T_e^*}\, \diff{t}
    + \alpha\, \lambda_2
    - \frac{1}{\reno} \lapl[x]{\lambda_2}
    - u^* \cdot \grad[x]{\lambda_2} &+{} \\
    \label{eq:ad_pde2}
    + \lambda_2 \cdot \grad[x]{u^*}
    - \grad[x]{\lambda_3}
    &= 0, \\
    \label{eq:ad_pde3}
    \diver[x]{\lambda_2} &= 0,
  \end{align}
  with initial and final conditions $\lambda_6(x) = \lambda_1(x, t_i)$ and $\lambda_1(x,t_i + T') = 0$ for each $x \in \Omega$, and boundary conditions $\lambda_1(x,t) = 0$ and $\lambda_2(x) = 0$ for each $x \in \bdry{\Omega}$ and $t \in [t_i, t_i + T']$.
\end{theorem}
The proof to Thm.~\ref{thm:estimate_lg1} is an extension of~\cite[Theorem 1.17]{ito2008lagrange}, and therefore omitted. 
Note that we also omitted the differential equations describing $\lambda_4$ and $\lambda_5$, since they are irrelevant to our results in this paper.
Using the result in Thm.~\ref{thm:estimate_lg1} we obtain closed-form formulas for the Fr\'echet derivatives of $J_c$ as defined in eqs.~\eqref{eq:mpc_objective}.
\begin{theorem}
  \label{thm:frechet}
  Let $J_c(T_e,u,p,g_{T_e},g_u)$ be differentiable.
  Then, the Fr\'echet derivatives of $J$ with respect to $g_{T_e}$ and $g_u$ exist.
  Moreover, for each $\delta g_{T_e} \in L^2\pb{[t_0,t_f];H^1(\Omega)}$, and $\delta g_u \in L^2(\Omega) \times L^2(\Omega)$:
  \begin{align}
    \hspace{-.7em}
    \label{eq:gte_deriv}
    &\dprodb{\fderiv[g_{T_e}]{J_c}}{\delta g_{T_e}}_{\Omega \times [t_0, t_f]}
    = \dprod{\grad[g_{T_e}]{J_c} - \lambda_1}{\delta g_{T_e}}_{\Omega \times [t_0, t_f]},\\
    \label{eq:gu_deriv}
    &\dprodb{\fderiv[g_{u}]{J_c}}{\delta g_{u}}_{\Omega}
    = \dprod{\grad[g_{u}]{J_c} - \lambda_2}{\delta g_u}_{\Omega}.
  \end{align}
\end{theorem}
The formulas above result from applying standard linear variations to the
objective functions, as shown in~\cite{he2016gradient} and we omitt
the detial derivatives in this paper.

\subsection{Gradient-Based Optimization}
\label{subsec:grad_method}
Using the closed-form formulas for the Fr\'echet derivatives of $J_e$
and $J_c$, as defined in~\cite{he2016gradient} and Thm.~\ref{thm:frechet}, we build a
gradient-based optimization algorithm to solve the problems in
eqs.~\eqref{eq:cost_estimate} and~\eqref{eq:cost_control} inspired in
a project-gradient approach~\cite[Ch.~18.6]{nocedal2006numerical}. 
Recall that projected gradient is a iterative method where a sequence
of quadratic programming (QP) problems are solved. 

Thus, we find descent directions $\delta g_{T_e}$ and $\delta g_u$ as the unique solutions of the following QP with value $V_c$:
\begin{equation}
  \label{eq:qp_grad_control}
  \begin{aligned}
    V_c \p{g_{T_e}, g_u} \!=\!\!
    \min_{\delta g_{T_e}, \, \delta g_u}
    &\!\!\!\dprod{\fderiv[g_{T_e}]{\!J_c}}{\delta g_{T_e}}_{\Omega \times [t_0,t_f]} \hspace{-.4em}+ \!\dprod{\fderiv[g_{u}]{J_c}}{\delta g_{u}}_{\Omega} \!+ \\
    &\ + M_{T_e}(\delta g_{T_e}, \delta g_{T_e}) + M_u(\delta g_u, \delta g_u),\\
    \text{s.t.:}\ \
    & (g_{T_e} + \delta g_{T_e})(x,t) \in \spb{\underline{g_{T_e}}, \overline{g_{T_e}}}, \\
    & (g_{u} + \delta g_{u})(x) \in \spb{\underline{g_{u}}, \overline{g_{u}}},\ \forall x\ \forall t,
  \end{aligned}
\end{equation}
where $M_{T_e}$ and $M_u$ are positive-definite bilinear operators.

Once we have found descent directions using the QP problems in eqs.~\eqref{eq:qp_grad_control}, we use Armijo's method~\cite{armijo1966minimization} to compute step sizes that guarantee the convergence of our algorithms to points satisfying first-order necessary optimality conditions.
That is, we find the largest scalar $\beta_c \in [0,1]$ such that, 
if $V_c(g_{T_e},g_u) < 0$ then
\begin{multline}
  \label{eq:armijo_control}
  J_c\p{g_{T_e} + \beta_c\, \delta g_{T_e}, g_{u} + \beta_c\, \delta g_{u}} - J_c\p{g_{T_e}, g_u} \leq\\
  \leq a \, \beta_c\, V_c\p{g_{T_e}, g_u},
\end{multline}
where $a \in (0,1)$ is a parameter.

\begin{figure}[tp]
  \centering
  {\small
    \begin{algorithmic}[1]
      \Require $i = 0$,~$t_0 = 0$,~$\Delta > 0$,~$T > \Delta$,~$T' > \Delta$,~$\epsilon_{\text{tol}} \geq 0$, $\pi_0 \in H^1(\Omega)$, $\theta \in \R^n$, $g_{T_e} \in L^2\pb{[t_0,t_f],H^1(\Omega)}$, $g_u \in L^2(\Omega) \times L^2(\Omega)$, $\setb{T_{e,k}^*|_{[-T,0]}}_{k=1}^{n_t}$.
      \Loop
      \State Update $\pi_0$ and $\theta$ by eqs.~\eqref{eq:cost_estimate} and~\cite[Algorithm 1]{he2016gradient}. \label{algo:estimation_after}
      \Loop \label{algo:control_begin}
      \Comment Solve control problem~\eqref{eq:cost_control} for $t \in [t_i,t_i+T']$.
      \State Compute $T_e$, $u$, and $p$ by solving~\eqref{eq:temp} to~\eqref{eq:bnd_u1}.
      \State Compute $\lambda_j$, $j \in \set{1,2,3,6}$, by solving~\eqref{eq:ad_pde1} to~\eqref{eq:ad_pde3}.
      \State Compute $\fderiv[g_{T_e}]{J_c}$ and $\fderiv[g_u]{J_c}$ in~\eqref{eq:gte_deriv} and~\eqref{eq:gu_deriv}.
      \State Compute $V_c(g_{T_e}, g_u)$, $\delta g_{T_e}$, and $\delta g_u$ by solving~\eqref{eq:qp_grad_control}.
      \If{$V_c(g_{T_e}, g_u) \geq -\epsilon_{\text{tol}}$}
      \State Go to line~\ref{algo:control_after}.
      \EndIf
      \State Compute $\beta_c$ in~\eqref{eq:armijo_control}.
      \State Update $g_{T_e} \leftarrow g_{T_e} + \beta_c\, \delta g_{T_e}$ and $g_u \leftarrow g_u + \beta_c\, \delta g_u$.
      \EndLoop \label{algo:control_end}
      \State Apply $g_{T_e}|_{[t_i,t_i+\Delta]}$ and $g_u|_{[t_i,t_i+\Delta]}$. \label{algo:control_after}
      \State Update $t_i \leftarrow t_i + \Delta$ and $i \leftarrow i + 1$.
      \State Store $\setb{T_{e,k}^*|_{[t_i-T,t_i]}}_{k=1}^{n_t}$.
      \EndLoop
    \end{algorithmic}
  }
  \caption{Model predictive HVAC estimation and control algorithm.}
  \label{fig:algo}
\end{figure}

We summarize our HVAC estimation and control algorithm in Fig.~\ref{fig:algo}.
Note that the parameter $\epsilon_{\text{tol}}$ defines the convergence tolerance of our optimization algorithm.
In practical applications it is recommended to choose $\epsilon_{\text{tol}} > 0$, but the convergence of the optimization algorithm is guaranteed only when $\epsilon_{\text{tol}} = 0$, as shown below.

\begin{theorem}
  \label{thm:algorithm_converge}
  Consider the algorithm in Fig.~\ref{fig:algo} and assume that $\epsilon_{\text{tol}} = 0$.
  Then the sequences generated by the sub-method in
  lines~\ref{algo:control_begin} to~\ref{algo:control_end} accumulates
  in points nullifying $V_c$. 
\end{theorem}

The proof of Theorem~\ref{thm:algorithm_converge} is a direct
extension of the result in~\cite[Ch.~4]{polak2012optimization} and we
omitt it in this paper.

\section{SIMULATED RESULTS}
\label{sec:results}

\begin{figure}[tp]
  \centering
  \includegraphics[width=.7\linewidth, trim=0 22 0 23, clip]{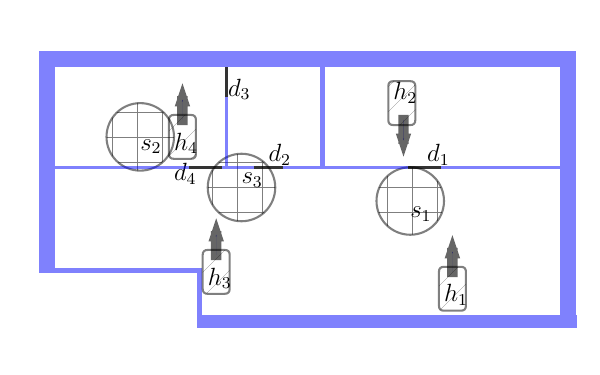}
  \caption{Simulated apartment floor plan.}
  \label{fig:floor}
\end{figure}

We applied our method to a simulated St.\,Louis area apartment, whose floor plan is shown in Fig.~\ref{fig:floor}.
The apartment has an area of $7.6 \times 16.8 \unit{m^2}$ ($1375 \unit{sq\,ft}$, approx.), has $n_d = 4$ doors, denoted $\set{d_i}_{i=1}^{n_d}$, $n_t = 3$ thermostats, labeled $\set{s_i}_{i=1}^3$, and is equipped with four HVAC vents, labeled $\set{h_i}_{i=1}^4$.
We assume that each vent works independently, endowed with a fan acting on an $1 \times 0.5 \unit{m^2}$ area.

Our CFD model is governed by the constants $\reno = 10^2$, $\kappa_0 = 10^{-2}$, $\alpha_w = 10^3$, $\kappa_w = 10^{-4}$, and $p_A = 101.3 \unit{kPa}$.
The sensor sensitivity is $r = 1.0 \unit{m}$, and the parameters in~\eqref{eq:mpc_objective} are $\eta'_0 = 0.1$, and $\eta'_1 = 0.15$.
The bilinear operators in~\eqref{eq:qp_grad_control} are estimated using the BFGS method~\cite[Ch.~6.1]{nocedal2006numerical}.
Our simulations are implemented using \emph{Python}, and the each PDE is solved using a FEM discretization computed with tools from the \emph{FEniCS} project~\cite{logg2012automated}.
The floor plan was discretized into $6276$ discrete elements.

\subsection{Variable door configuration and target location}
\label{subsec:simulation_1}
In this experiment we study the accuracy of our estimation method discussed by eqs.~\eqref{eq:cost_estimate} under changing door configurations, while our control methods maximizes efficient and comfort in 12~different target areas uniformly distributed across the apartment.

We compare the performance of our method against a scheme consisting of a similar controller but no door configuration estimator, in order to validate the improvement to controller by using our coupled estimator.
The comparison controller assumes that the doors are always closed.
The atmospheric temperature is set to $T_A = 5 \unit{^{\circ} C}$, the clothing insulation to $0.155 \unit{^{\circ}C\,m^2/W}$, and the metabolic rate to $64.0 \unit{W/m^2}$, resulting in an outdoor PMV index of $-4.1$ (i.e., very cold according to~\cite{ashrae2010}).
The simulation begins with all four doors closed, and at time $50 \unit{s}$ all the doors are simultaneously open.

\begin{figure*}[tp]
  \centering
  \subfloat[%
    PMV distribution for specific target area.%
  ]{%
    \label{fig:PMV_exp1_1}%
    \includegraphics[width=.31\linewidth,trim=32 0 6 0,clip]{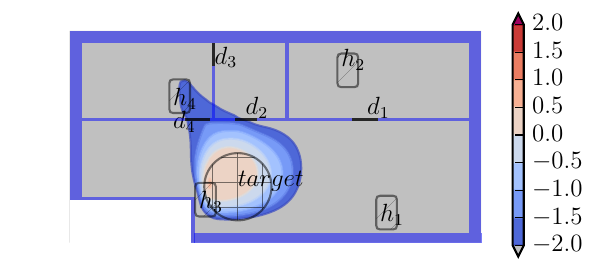}%
  }%
  \hfill%
  \subfloat[%
    PMV distribution for specific target area.%
  ]{%
    \label{fig:PMV_exp1_2}%
    \includegraphics[width=.31\linewidth,trim=32 0 6 0,clip]{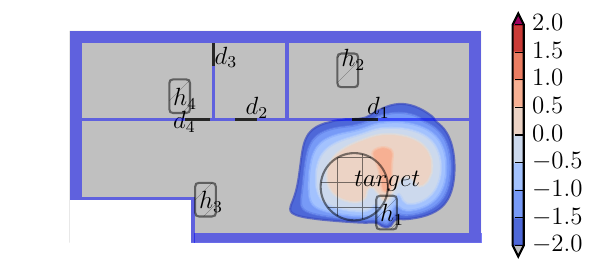}%
  }%
  \hfill%
  \subfloat[%
    Temperature error distribution when no estimator is used for specific target area.%
  ]{%
    \label{fig:tem_err}%
    \includegraphics[width=.31\linewidth,trim=32 0 6 0,clip]{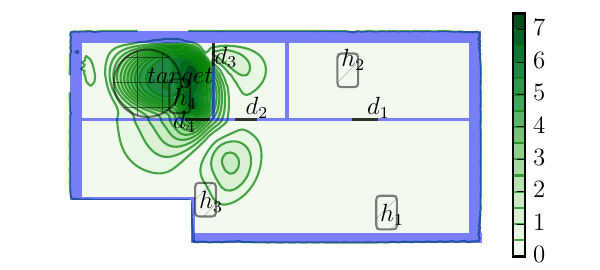}%
  }%
  \caption{%
    Results of representative experiments in Sec.~\ref{subsec:simulation_1}.
    Each plot was generated at $t=900\unit{s}$ using a different target area.
  }
  \label{fig:simulation1_distribution}
\end{figure*}

Figs.~\ref{fig:PMV_exp1_1} and~\ref{fig:PMV_exp1_2} show how our method adapts to changing target areas when the door estimator is enabled, quickly reaching comfortable conditions after a few minutes while reducing energy consumption in unused areas.

\begin{figure}[tp]
  \centering
  \subfloat[%
    Average absolute PMV value in target area.
    Left~(M): using our estimator.
    Right~(F): assuming closed doors.
  ]{%
    \label{fig:PMV_exp1}%
    \hspace{.085\linewidth}%
    \includegraphics[width = 0.3\linewidth,trim=0 20 10 10,clip]{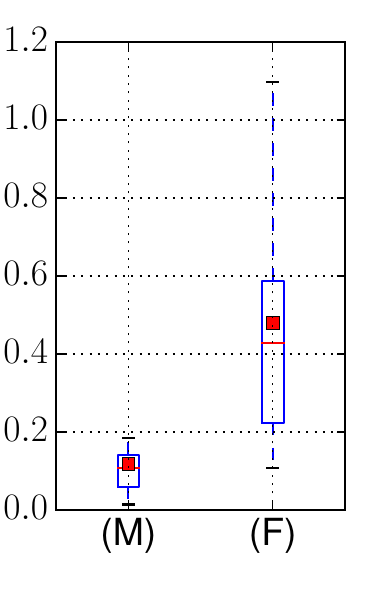}%
    \hspace{.085\linewidth}%
  }%
  \hfill%
  \subfloat[%
    Average estimation error over all target areas, $\frac{1}{n_d} \sum_{i=1}^{n_d} \absb{\theta_i - \hat{\theta}_i}$, where $\hat{\theta}$ is our estimator.]{%
    \label{fig:theta_exp1}%
    \includegraphics[width = 0.47\linewidth,trim=5 0 10 10,clip]{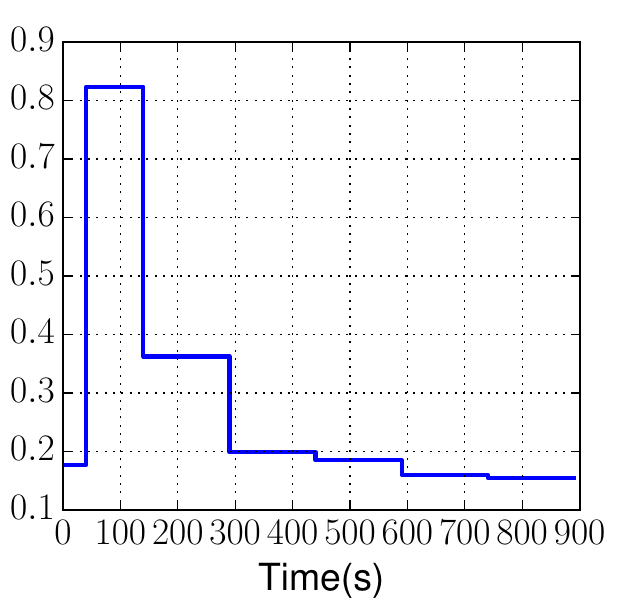}%
  }%
  \caption{%
    Results of the simulations in Sec.~\ref{subsec:simulation_1} using 12 different target areas.
  }
  \label{fig:simulation1_stats}
\end{figure}

Fig.~\ref{fig:tem_err} shows the impact of our door configuration
estimation method in the temperature distribution.
Indeed, the changes in airflow due to different door configurations
are significant enough to produce a sizable temperate estimation
error, which in turns leads to uncomfortable conditions as shown in
Fig.~\ref{fig:PMV_exp1}.
Fig.~\ref{fig:theta_exp1} shows that, after the change in door
configuration at $t=50\unit{s}$, it takes our door configuration
estimator about $250\unit{s}$ to recalculate the new door
configuration, which is compatible with the time constant of the
system's dynamics.

\subsection{Simulation with different personal variables related to PMV}
\label{subsec:simulation_2}

\begin{figure}[tp]
  \centering
  \subfloat[%
    Energy usage (\unitns{kW\, h}) vs.\ meta\-bolic rate (\unitns{W/m^2}) with clothing insulation $0.155~\unit{^{\circ}C\,m^2/W}$.
  ]{%
    \label{fig:energy_meta}%
    \includegraphics[width = 0.47\linewidth,trim= 15 5 10 10,clip]{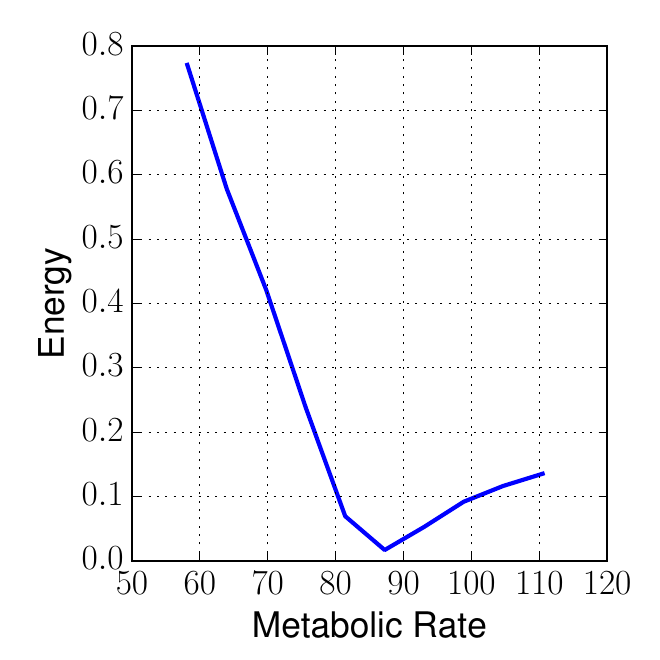}%
  }%
  \hfill%
  \subfloat[%
    Energy usage (\unitns{kW\, h}) vs.\ clothing insulation (\unitns{^{\circ}C\,m^2/W}) with metabolic rate $81 \unit{W/m^2}$.
  ]{%
    \label{fig:energy_clo}%
    \includegraphics[width = 0.47\linewidth,trim= 15 5 10 10,clip]{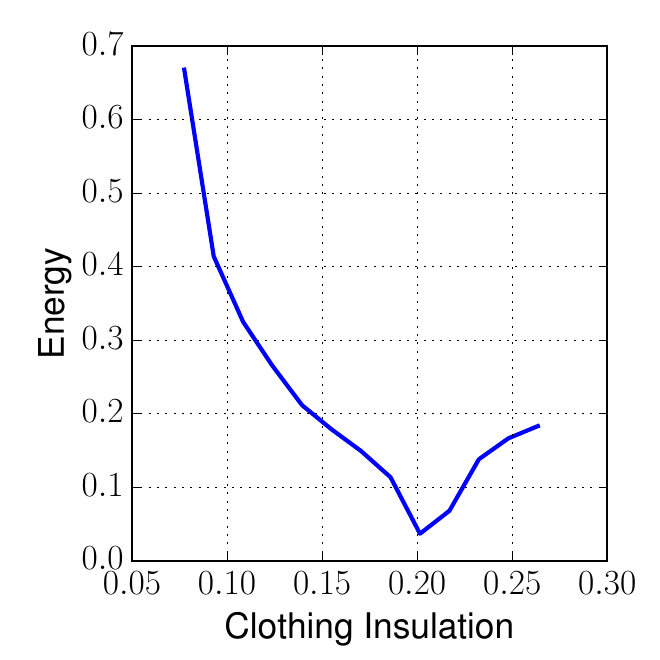}%
  }%
  \caption{%
    Results of the simulations in Sec.~\ref{subsec:simulation_2}.
  }
  \label{fig:exp_stats}
\end{figure}


In this simulation we show how our climate control scheme is capable of not only adapting to changes in building configuration, but also to changes in the occupancy patterns and behaviors.
In particular, we consider the cases where metabolic rate and clothing insulation are known variables, and feed that information back to the controller.
While this information is not typically known in a traditional climate control setup, it is possible to obtain it via wearable sensors or smartphone surveys~\cite{jin2014environment, cleland2013optimal, yang2010review}.

Figs.~\ref{fig:energy_meta} and~\ref{fig:energy_clo} show the impact that changes in metabolic rate and clothing insulation have in the energy consumption when using our climate control scheme, respectively.
We ran simulations over a period of $600 \unit{s}$, with atmospheric temperature set to $20 \unit{^{\circ} C}$.
Our choice of atmospheric temperature is such that the PMV value is negative when occupants have low metabolic rate and clothing index, and is positive when occupants have high metabolic rate and clothing index.
Thus, as we vary the values for metabolic rate and clothing insulation, the HVAC system automatically switches from heating to cooling.
As shown in both figures, both variables have a significant impact in energy consumption, particularly when heating is required due to low metabolic rates or clothing insulation.
Therefore, our scheme not only results in a more comfortable
environment under different occupancy patterns, it also opens the door
to the semi-autonomous operation of a building where the climate
control system suggests, in real-time, changes in clothing insulation
or occupant behavior that would help saving energy.

\subsection{Memory and time usage analysis}
\label{subsec:simulation_3}

\begin{figure}[tp]
  \centering
  \subfloat[%
    Memory usage (\unitns{GB}) vs.\ number of FEM elements. %
  ]{%
    \label{fig:memory_element}%
    \hspace{.02\linewidth}%
    \includegraphics[width=.43\linewidth, trim=5 5 15 15, clip]{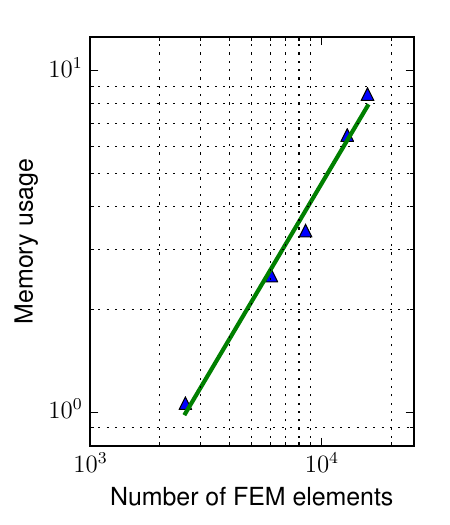}%
    \hspace{.02\linewidth}%
  }%
  \hfill%
  \subfloat[%
    Computation time (\unitns{sec}) vs.\ number of FEM elements.%
  ]{%
    \label{fig:time_element}%
    \hspace{.02\linewidth}%
    \includegraphics[width=.43\linewidth, trim=5 5 15 15, clip]{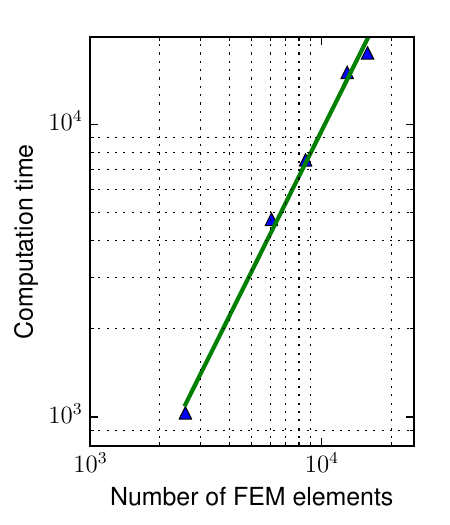}%
    \hspace{.02\linewidth}%
  }%
  \caption{%
    Results of the experiments in Sec.~\ref{subsec:simulation_3}.
  }
  \label{fig:simulation3}
\end{figure}


Fig.~\ref{fig:memory_element} shows the maximum memory usage under different number of FEM elements, $n_e$, resulting in a memory complexity roughly proportional to $n_e^{1.14}$.
Similarly, Fig.~\ref{fig:time_element} shows the average CPU time of each iteration of the algorithm in Fig.~\ref{fig:algo} under different number of FEM elements, resulting in a time complexity roughly proportional to $n_e^{1.59}$.

\begin{figure}
  \centering
  \includegraphics[width=.6\linewidth, trim=20 18 15 10, clip]{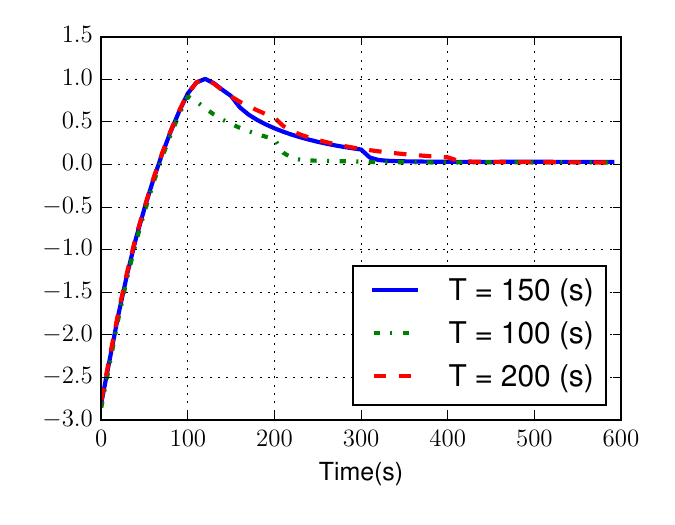}%
  \caption{%
    Average PMV value in target area for different control horizons.
  }
  \label{fig:PMV_time}%
\end{figure}

While our scheme is still not fast enough to run in real-time, our CPU time is within less than an order of magnitude of the desired speed.
Also, real-time implementations are a function of the estimation and control horizons, denoted $T$ and $T'$ in Fig.~\ref{fig:algo}.
As shown in Fig.~\ref{fig:PMV_time}, there is not a large impact when
longer horizons are used, mostly due to the slow dynamic that climate
variables have in general.

\section{CONCLUSIONS}
In this paper, we introduce a comfort aware building climate control using distributed-parameter CFD model.
Beside temperature, the target function of the controller considers multi-variables.
Our controller coupled with our estimator show their accuracy and energy-efficiency to estimate and improve indoor climate.
Future study will focus on accelerating algorithms to solve large-scale nonlinear optimal control problem to improve its computation efficiency.

\bibliographystyle{IEEEtran}
\bibliography{bibfile}

\appendices
\section{Predicted Mean Vote Formula}
\label{appendix:pmv}

According to Fanger~\cite{fanger1970thermal} and
ASHRAE~\cite{ashrae2010}, the PMV index is computed by
\begin{multline}
  \label{eq:pmv}
  \pmv = \ \p{0.303\, \exp\p{-0.036\, M} + 0.028}\, \pB{\pn{M - W} +{} \\
    - 3.05\cdot 10^{-3}\, \pb{5733 - 6.99\, \pn{M - W} - p_a} +{} \\
    - 0.42\, \pn{M - W - 58.15} - 1.7\cdot 10^{-5}\, M\, \pn{5867 - p_a} +{} \\
    - 0.0014\, M\, \paren{34 - T_e} - f_{cl}\, h_c\, \paren{T_{cl} - T_e} +{} \\
    - 3.96\cdot 10^{-8}\, f_{cl}\, \pb{\pn{T_{cl} + 273}^4 - \pn{T_r + 273}^4}}.
\end{multline}
Where $M$ and $W$ are the metabolic rate and external work, both in
$\unit{W/m^2}$.
The external work normally is around zero~\cite{awbi2003ventilation,
  ashrae2010}, and the average value of human's sedentary activity and
standing activity is $70 \unit{W/m^2}$ and $93 \unit{W/m^2}$
respectively.
The occupant's metabolic rate can be obtained either from their
wearable devices or some posterior
estimation~\cite{jin2014environment, cleland2013optimal,
  yang2010review}.

The term $p_a$ corresponds to the partial water vapor pressure,
measured in Pascals.
According to~\cite{Awbi1986domestic, sasic2004ham}, the specific
humidity inside house, $w_i$ with unit $\unit{kg \cdot kg^{-1}}$, can
be expressed by
\begin{equation}
  \label{eq:moisture_steady}
  w_i = \frac{w_o \rho A\, g_u + m_g}{\rho A\, g_u}.
\end{equation}
Where $w_{o}$ is the specific humidity comes out of the HVAC,
$\rho$ is the inside air's density. 
$m_g$ is the rate of moisture generation within the building with unit
$\unit{kg \cdot s^{-1}}$,
$A$ is the size of the fan with unit $\unit{m^2}$.
In the computation, we ignore the moisture diffusion through the fabric
material~\cite{awbi2003ventilation}. 
Under ideal air condition, $p_a = 1.608\, p_o\, w_i$, where the mixed
air's pressure $p_o$ is the standard atmosphere, $1.013 \cdot 10^5
\unit{Pa}$.

$T_e$ and $T_r$ are the air temperature and mean radiant temperature,
both with unit $\unit{^{\circ}C}$. 
Experiments results~\cite{Walikewitz2015difference} show that their
indoor distributions are close. 
In order to simplify the PMV index for optimization computation and
noting that most buildings typically do not have sensors to
continually measure the mean radiant temperature, we set the $T_r$
equal to the air temperature, $T_e$. 
Moreover, since the range of temperatures in the indoor environment is
small, the difference of the fourth power terms can be adequately
replaced by a lower-order difference~\cite{awbi2003ventilation}:
\begin{multline}
  \label{eq:radiation}
    3.96 \cdot 10^{-8}\, f_{cl}\, \pb{\p{T_{cl} + 273}^4 - \paren{T_r + 273}^4}
    \approx \\
    \approx 4.6\, f_{cl}\, \p{1 + 0.01 T_r}\, \p{T_{cl} - T_r}.
\end{multline}
The clothing surface temperature is approximated
as~\cite{awbi2003ventilation, fanger1970thermal}: 
\begin{multline}
  \label{eq:t_cl}
  T_{cl} =
  - 0.155\, I_{cl}\, f_{cl}\, \pB{4.6\, \p{1 + 0.01 T_r}\, \p{T_{cl} - T_r} +{} \\
    + h_c\, \p{T_{cl} - T_e}} - 0.028\, \p{M - W} + 35.7.
\end{multline}
Where we approximate the radiation term with
equation~\eqref{eq:radiation}, then we can derive an explicit formula
of $T_{cl}$ and get rid of the iteration numerical solving process. 
$h_c$ is the convective heat transfer
coefficient~\cite{fanger1970thermal, colin1967experimental} and 
is approximated as the natural convective heat transfer coefficient,
$h_{cn}$~\cite{cigler2012optimization, federspiel1992}. 
The parameter $f_{cl}$ is equal to $1.0 + 1.29\, I_{cl}$ when $I_{cl}
\leq 0.078$, otherwise $1.05 + 0.645\, I_{cl}$, where $I_{cl}$ is the
clothing insulation index, in~$\unit{m^2\, ^{\circ}C / W}$.

After the simplications and derivatives above, the PMV index follows
equation~\eqref{eq:pmv_simplification} where the parameters,
$\{a_i\}_{i=0}^{10}$, are
\begin{equation}
  \begin{aligned}
    & a_0 = \p{0.303\, \exp\p{-0.036\, M} + 0.028}\,
    \pB{\pn{M - W} +{} \\
      & \hspace{1.0em}
      - 3.05\cdot 10^{-3}\, \pb{5733 - 6.99\, \pn{M - W}}
      - 0.0476\, M
      +{} \\
      & \hspace{1.0em}
      - 0.42\, \pn{M - W - 58.15} - 9.974 \cdot 10^{-2}\, M +{} }, \\
    & a_1 = \p{0.303\, \exp\p{-0.036\, M} + 0.028} \cdot{} \\
    & \hspace{1.0em}
    \cdot \p{3.05\cdot
      10^{-3} + 1.7 \cdot 10^{-5} M} w_o \rho A, \\
    & a_2 = \p{0.303\, \exp\p{-0.036\, M} + 0.028} \cdot{} \\
    & \hspace{1.0em}
    \cdot \p{3.05\cdot
      10^{-3} + 1.7 \cdot 10^{-5} M} m_g, \\
    & a_3 = \rho A, \\
    & a_4 = \p{0.303\, \exp\p{-0.036\, M} + 0.028} \cdot{} \\
      &\hspace{1.0em}
    \cdot \, \p{.0014\, M + f_{cl}\, h_c + 4.6 f_{cl}}, \\
    & a_5 = 0.046 f_{cl} \p{0.303\, \exp\p{-0.036\, M} + 0.028}, \\
    & a_6 = \p{0.303\, \exp\p{-0.036\, M} + 0.028} \cdot{} \\
    & \hspace{1.0em}
    \cdot (4.6 + h_c)f_{cl}\, (-0.028 (M - W) + 35.7), \\
    & a_7 = \p{0.303\, \exp\p{-0.036\, M} + 0.028} \cdot{} \\
    & \hspace{1.0em}
    \cdot{} I_{c} f^2_{cl} \p{0.155h^2_c + 1.426 h_c + 3.45}, \\
    & a_8 = 1, \\
    & a_9 = I_c f_{cl} \p(0.713 + h_c), \text{ and} \\
    & a_{10} = 0.007 I_c f_{cl}.
  \end{aligned}
  \label{eq:pmv_as}
\end{equation}

\section{Existence and Uniqueness of Solutions to the CFD Model}
\label{appendix:existence}

In this section, we are going to show the solution's existence and uniqueness to the CFD model containing equations~\eqref{eq:temp}, \eqref{eq:init_t}, \eqref{eq:ns1}, \eqref{eq:ns2}, \eqref{eq:bnd_t} and~\eqref{eq:bnd_u1} in Sec.~\ref{subsec:cfd_model}.

First, let us describe our notation.
The definitions and notations of $L^2(\Omega)$ and the Sobolev space $H^1(\Omega)$ follow~\cite{girault2012finite}.
Let the subspace for $H^1(\Omega)$ with trivial boundary be denoted:
\begin{equation}
  H^1_0(\Omega) = \set{v \in H^1(\Omega) \mid v|_{\partial \Omega} = 0}.
\end{equation}
Let the divergence-free subspace of $H_0^1(\Omega) \times H_0^1(\Omega)$ be defined by:
\begin{equation}
  V_0 = \set{v \in H^1_0(\Omega) \times H^1_0(\Omega) \mid \grad[x]{} \cdot v = 0}.
\end{equation}
Let $H_0$ be the completion of ${V}_0$ w.r.t. $L^2$ norm, and is given by:
\begin{equation}
  H_0 = \set{v \in L^2(\Omega) \times L^2(\Omega) \mid \diver[x]{v} = 0,\ \text{and}\ v|_{\bdry{\Omega}} = 0}.
\end{equation}
The space $H_0$ is equipped with the $L^2$-norm, denoted $\norm{\cdot}_0$.
Let $V_1 = \set{\eta \in H^1(\Omega) \mid \eta|_{\partial \Omega} = 0}$, endowed with the norm $\norm{f}_1 = \norm{\grad[x]{f}}_0$ for each $f \in V_1$.
Also, let:
\begin{multline}
  L^2([t_0, t_0 + T]; V_1)
  = \setbb{u\colon [t_0, t_0 + T] \to V_1 \mid\\
    u\ \text{is measurable and}\ \int_{t_0}^{\mathrlap{t_0 + T}} \norm{u(t)}_1^2\, \diff{t} < \infty},
\end{multline}
and:
\begin{multline}
  C^0\pb{[t_0, t_0 + T]; L^2(\Omega)}
  = \setb{u\colon [t_0, t_0 + T] \to L^2(\Omega) \mid\\
    u(t)\ \text{is measurable and continuous for a.e.}\ t}.
\end{multline}


Now we focus our attention on the existence and uniqueness of the CFD system in eqs.~\eqref{eq:temp} to~\eqref{eq:bnd_u1}, stated in Thms.~\ref{thm:wk_existence1} and~\ref{thm:wk_unique1}.
In order to prove these theorems, we begin studying the weak solution PDE subsystem in eqs.~\eqref{eq:ns1}, \eqref{eq:ns2} and boundary condition~\eqref{eq:bnd_u1}.
The techniques used in our proofs follow closely those previously published results~\cite{girault2012finite, desai1994optimal, ito1998optimal}.

According to the Hopf extension~\cite[Lemma I.4.2.3]{girault2012finite}, for each $\epsilon > 0$ there exists a function $\bar{u} \in H^1(\Omega) \times H^1(\Omega)$ such that $\diver[x]{\bar{u}} = 0$, $\bar{u}|_{\bdry{\Omega}} = 0$, and $\absb{\dprodb{\p{v \cdot \grad[x]}\, \bar{u}}{v}_{\Omega}} \leq \epsilon \abs{v}_1^2$ for each $v \in V_0$.
Then, any function $u \in H^1(\Omega) \times H^1(\Omega)$ satisfying
the boundary condition~\eqref{eq:bnd_u1} can be written as $u = w +
\bar{u}$, for some $w \in V_0$.
For simplification, define $\Sigma$ as $H^1(\Omega) \times H^1(\Omega)$.

We use a fixed-point theorem to show the existence of weak solutions of eq.~\eqref{eq:ns1}.
Consider the following weak formulation of~\eqref{eq:ns1} for each $\phi \in V_0$:
\begin{multline}
  \label{eq:wk_ln1}
    \dprod{\alpha(x) w}{\phi}_{\Omega}
    + \frac{1}{\reno} \dprod{\grad[x]{u}}{\grad[x]{\phi}}_{\Omega}
    + \dprod{\hat{u} \cdot \grad[x]{w}}{\phi}_{\Omega}
    + \\
    + \dprod{u \cdot \grad[x]{\bar{u}}}{\phi}_{\Omega}
    = \dprod{g_u}{\phi}_{\Omega}.
\end{multline}
Consider the mapping $\mathcal{S}$ such that given $\hat{u}$, it returns $u$, the solution to the equation~\eqref{eq:wk_ln1}, as $\mathcal{S}(\hat{u}) = u$.
Then, the fixed point of the mapping $\mathcal{S}$ is a weak solution to the equation~\eqref{eq:ns1}.
Let us define the following bilinear operator, for each $w, \phi \in V_0$:
\begin{multline}
  \sigma_0(w, \phi)
  = \dprod{\alpha(x) w}{\phi}_{\Omega}
  + \frac{1}{\reno}\, \dprod{\grad[x]{u}}{\grad[x]{\phi}}_{\Omega} +\\
  + \dprod{\hat{u} \cdot \grad[x]{w}}{\phi}_{\Omega}
  + \dprod{u \cdot \grad[x]{\bar{u}}}{\phi}_{\Omega}.
\end{multline}

The proof of Lemmas~\ref{lem:sigma0_bound} to~\ref{lem:su_compact} below can be found in~\cite{desai1994optimal}.
\begin{lemma}
  \label{lem:sigma0_bound}
  $\sigma_0$ is bounded.
\end{lemma}

\begin{definition}
  A bilinear form $a\colon V \times V \to \R$ is called \emph{coercive} if there exists a constant $C > 0$ such that for each $x \in V$, $\abs{a(x,x)} \geq C \norm{x}_V^2$.
\end{definition}

\begin{lemma}
  \label{lem:sigma0_coercive}
  $\sigma_0$ is coercive.
\end{lemma}

\begin{lemma}
  \label{thm:wk_ln1}
  For each $\phi \in V_0$, the equation~\eqref{eq:wk_ln1} has a unique weak solution $u$.
  Moreover, there exists $w \in V_0$ such that $u$ can be separated as $u = w + \bar{u}$.
\end{lemma}

\begin{lemma}
  \label{lem:su_bound}
  For each $\hat{u} \in V_0$, the solution $u = \mathcal{S}(\hat{u})$
  to equation~\eqref{eq:wk_ln1} is always bounded.
\end{lemma}

\begin{lemma}
  \label{lem:su_compact}
  The mapping ${\cal S} \colon \Sigma \to \Sigma$ is compact.
\end{lemma}

Based on these lemmas, we can derive the existence of weak solutions to eq.~\eqref{eq:ns1}.
\begin{lemma}
  \label{thm:wku_existence}
  There exists at least one fixed point, say $u$, for the mapping $\mathcal{S}$, such that $u \in V_0 + \bar{u}$.
\end{lemma}
The proof of Lemma~\ref{thm:wku_existence} follows directly from Lemma~\ref{lem:su_compact}, as shown in~\cite[Theorem 1.J]{troianiello2013elliptic}.

We now turn our attention to the existence and uniqueness of weak solutions of eq.~\eqref{eq:temp}.
Let us define bilinear operator $\sigma_1$ for each $t \in [t_0, t_0 + T]$ by:
\begin{equation}
  \sigma_1(T_e, \xi)
  = \dprod{\kappa(x)\, \grad[x]{T_e}}{\grad[x]{\xi}}_{\Omega}
  + \dprod{u \cdot \grad[x]{T_e}}{\xi}_{\Omega}.
\end{equation}
Then, given $t \in [t_0, t_0 + T]$, the weak formulation of eq.~\eqref{eq:temp} can be written as follows, where $\xi \in V_1$:
\begin{equation}
  \dprod{\deriv{T_e}{t}}{\xi}_{\Omega}
  + \sigma_1(T_e, \xi)
  = \dprod{g_{T_e}}{\xi}_{\Omega}.
\end{equation}

\begin{lemma}
  \label{lem:sigma1_bound}
  $\sigma_1$ is bounded.
\end{lemma}
\begin{IEEEproof}
  Since $\kappa(x)$ and $u(x)$ are bounded, using H\"older's inequality we can show that $\sigma_1$ is bounded using the same argument as in Lemma~\ref{lem:sigma0_bound}.
\end{IEEEproof}

\begin{lemma}
  \label{lem:sigma1_coercive}
  $\sigma_1$ is coercive.
\end{lemma}
\begin{IEEEproof}
  First, according to equation~\eqref{eq:ns2} and derivative by part,
  $\dprod{u \cdot \grad[x]{T_e}}{T_e}_{\Omega} = 0$.
  Then by Poincare-Friedrichs inequality~\cite{vohralik2005discrete}, there
  exists $C > 0$ such that: 
  \begin{equation}
    \begin{aligned}
      \abs{\sigma_1(T_e, T_e)}
      &= \absb{\dprod{\kappa(x)\, \grad[x]{T_e}}{\grad[x]{T_e}}_{\Omega}
        + \dprod{u \cdot \grad[x]{T_e}}{T_e}_{\Omega}}\\
      &= \absb{\dprod{\kappa(x)\, \grad[x]{T_e}}{\grad[x]{T_e}}_{\Omega}}\\
      &\geq C\, \norm{T_e}_1^2,
    \end{aligned}
  \end{equation}
  as desired.
\end{IEEEproof}

\begin{IEEEproof}[Proof of Thm.~\ref{thm:wk_existence1}]
  The existence of weak solutions to eq.~\eqref{eq:temp} follows by Lemmas~\ref{lem:sigma1_bound} and~\ref{lem:sigma1_coercive}, together with~\cite[Theorem 11.1.1]{quarteroni2008numerical}.
  The existence of weak solutions to eqs.~\eqref{eq:ns1} and~\eqref{eq:ns2} follows by Lemma~\ref{thm:wku_existence}.
\end{IEEEproof}

Finally, we prove a sufficient condition for the uniqueness of the weak solutions to our CFD system.
\begin{IEEEproof}[Proof of Thm.~\ref{thm:wk_unique1}]
  Note that the uniqueness of the weak solution to eq.~\eqref{eq:temp} is guaranteed together with the existence result in~\cite[Theorem 11.1.1]{quarteroni2008numerical}.

  Suppose $u_1$ and $u_2$ are two different weak solutions to eq.~\eqref{eq:ns1} with the same boundary conditions.
  Let $\tilde{u} = u_1 - u_2$, then for each $\phi \in V_0$:
  \begin{multline}
    \dprod{\alpha(x)\, \tilde{u}}{\phi}_{\Omega}
    + \frac{1}{\reno}\, \dprod{\grad[x]{\tilde{u}}}{\grad[x]{\phi}}_{\Omega}
    + \dprod{u_1 \cdot \grad[x]{\tilde{u}}}{\phi}_{\Omega} +\\
    + \dprod{\tilde{u} \cdot \grad[x]{u_2}}{\phi}_{\Omega}
    = 0.
  \end{multline}
  Let us set $\phi = \tilde{u}$, then:
  \begin{multline}
    \dprod{\alpha(x)\, \tilde{u}}{\tilde{u}}_{\Omega}
    + \frac{1}{\reno}\, \dprod{\grad[x]{\tilde{u}}}{\grad[x]{\tilde{u}}}_{\Omega}
    + \dprod{u_1 \cdot \grad[x]{\tilde{u}}}{\tilde{u}}_{\Omega} +\\
    + \dprod{\tilde{u} \cdot \grad[x]{u_2}}{\tilde{u}}_{\Omega}
    = 0.
  \end{multline}
  Thus we can bound the norm of $\tilde{u}$ as:
  \begin{equation}
    \label{eq:tilde_u}
    \begin{aligned}
      \frac{1}{\reno} \dprod{\grad[x]{\tilde{u}}}{\grad[x]{\tilde{u}}}_{\Omega}
      &= - \dprod{\alpha(x)\, \tilde{u}}{\tilde{u}}_{\Omega} - \dprod{\tilde{u} \cdot \grad[x]{u_2}}{\tilde{u}}_{\Omega}\\
      &\leq \absb{\dprod{\tilde{u} \cdot \grad[x]{u_2}}{\tilde{u}}_{\Omega}}\\
      &\leq C\, \norm{u_2}_1\, \norm{\tilde{u}}_1^2,
    \end{aligned}
  \end{equation}
  where $C = \frac{\abs{\Omega}^{1/2}}{2}$, as shown in~\cite[Lemma 9.1.2]{galdi2011introduction}.

  Eq.~\eqref{eq:tilde_u} implies that $\p{\frac{1}{\reno} - C\, \norm{u_2}_1}\, \norm{\tilde{u}}_1^2 \leq 0$.
  Thus, if $\frac{1}{\reno} - C\, \norm{u_2}_1 \geq 0$ then $\norm{\tilde{u}}_1 = 0$, as desired.
\end{IEEEproof}


\end{document}